# The Cohomology of the McLaughlin Group and Some Associated Groups

ALEJANDRO ADEM†
Department of Mathematics
University of Wisconsin
Madison, WI 53706

R. JAMES MILGRAM†
Department of Math and Statistics
University of New Mexico
Albuquerque, NM 87131

## §0 Introduction

The group of order $2^7$, $Syl_2(M_{22}) = H$ seems quite remarkable – being the Sylow 2-subgroup of three of the sporadics, $M_{22}$, $M_{23}$, $McL$, as well as the simple group of Lie type $PSU_4(3)$. Also, it is a particularly nice index two subgroup of $Syl_2(Ly)$. We are interested in determining the mod(2) cohomology of these groups, and throughout this note all cohomology groups will be taken with $\mathbb{F}_2$ coefficients, so the coefficients will generally be suppressed. We use the standard Atlas (see [Co]) notation.

We have already determined the rings $H^*(M_{22}) \subset H^*(H)$, [AM2] and $H^*(M_{23}) \subset H^*(H)$, [M]. In this article we determine the rings

$$H^*(PSU_4(3)) \subset H^*(H) \text{ and } H^*(McL) \subset H^*(H).$$

Additionally, the group extension

$$N = 3 \cdot McL{:}2 \subset Ly$$

is a maximal subgroup of odd index in $Ly$, and we are able to determine $H^*(N)$ which contains $H^*(Ly)$ as a subring. This is very close to determining $H^*(Ly)$ itself. Indeed, what is required is a close study of the double coset decomposition of $Ly$ determined by $N$, which apparently is already reasonably well understood [Co].

The results are

THEOREM A: *There is a long exact sequence*

$$0 \longrightarrow \mathbb{F}_2[d_8, d_{12}](a_2, a_7, a_{11}, a_{14}) \longrightarrow H^*(PSU_4(3)) \longrightarrow H^*(2^4)^{\mathcal{A}_6} \oplus H^*(2^4)^{\mathcal{A}_6} \longrightarrow \mathbb{F}_2[d_8, d_{12}](1, \bar{w}_3, b_{15}, \bar{w}_3 b_{15}) \longrightarrow 0.$$

∎

† Both authors were partially supported by grants from the National Science Foundation and the ETH-Zürich. Also, the first author was partially supported by an NSF Young Investigator Award.

Here the subscript on an element denotes its dimension, $\mathbb{F}_2[x, \ldots]$ is the polynomial algebra on the stated generators, and the elements in parenthesis following the polynomial algebra are module generators. Also, in [AM4] the algebra of $\mathcal{A}_6$-invariants was determined as

$$\mathbb{F}_2[\bar{w}_3, \gamma_5, d_8, d_{12}](1, \gamma_9, b_{15}, \gamma_9 b_{15}).$$

In the case of $McL$ the result takes the form

THEOREM B: *There is a long exact sequence*

$$0 \longrightarrow \mathbb{F}_2[d_8, d_{12}](a_7, a_{11}) \longrightarrow H^*(McL) \longrightarrow \\ H^*(2^4)^{\mathcal{A}_7} \oplus H^*(2^4)^{\mathcal{A}_7} \longrightarrow \mathbb{F}_2[d_8, d_{12}](1, x_{18}) \longrightarrow 0.$$

∎

Again, in [AM4] the algebra of $\mathcal{A}_7$-invariants was determined as

$$\mathbb{F}_2[d_8, d_{12}, d_{14}, d_{15}](1, x_{18}, x_{20}, x_{21}, x_{24}, x_{25}, x_{27}, x_{45})$$

where $\mathbb{F}_2[d_8, d_{12}, d_{14}, d_{15}]$ is the Dickson algebra of $GL_4(2)$ invariants in $\mathbb{F}_2[x_1, y_1, z_1, w_1]$. In particular, $H^*(McL; \mathbb{F}_2)$ is 6 connected, and this is the largest connectivity of the 2-summand that we have seen yet in dealing with the sporadic groups.

Finally we have

THEOREM C: *The* mod(2) *cohomology of the maximal subgroup $N \subset Ly$ has the form*

$$\mathbb{F}_2[d_8, d_{12}](a_7, a_{11}) \oplus H^*(2^4)^{\mathcal{A}_7} \oplus \mathbb{F}_2[d_8, d_{12}, e_1](1, a_7, a_{11}, b_{18})e_1$$

*where $e_1$ is the generator coming from the split $\mathbb{Z}/2$-subgroup of $N$, and the remaining classes restrict to $H^*(McL)$.*

We also show in §6 that the first two summands in the display above are in the image of restriction from $H^*(Ly)$, but, as indicated, more work is necessary to complete the determination of $H^*(Ly)$. S. Smith and K. Umland are currently doing the necessary analysis of $Ly$, the results will appear in a sequel.

The techniques used here are quite a bit different from the methods usually used in the cohomology of finite groups. It is first critical to determine $H^*(PSU_4(3))$. To do this we use our knowledge of $H^*(SU_4(3))$ and a spectral sequence with $E_2$-term $\mathrm{Ext}_{H^*(\mathbb{Z}/4)}(H^*(SU_4(3)); \mathbb{F}_2)$ which converges to $H^*(PSU_4(3))$. This $E_2$-term is calculated and the $d_2$-differentials are calculated explicitly in §3. There may well be higher differentials in the spectral sequence, but $E_3 = E_\infty$ through at least dimension 8. This is enough to determine that there are two distinct classes in dimension 5 in $H^*(PSU_4(3))$. In turn, the existence of these two classes forces the remaining classes to have one and only one form – thus determining $H^*(PSU_4(3))$. Moreover, once this ring is determined, the calculation of $H^*(McL)$ is quite direct, as is the determination of $H^*(N)$.

The description of $H^*(PSU_4(3))$ is given in 3.7, the description of $H^*(McL)$ is given in 5.2 and 5.3, and the description of $H^*(N)$ is given in 6.1. The discussion in §4 gives the structure of the Bockstein spectral sequence for $H^*(PSU_4(3))$ which gives considerable information about the connections between the classes in these cohomology rings as well as determining the structure of the higher two–torsion there.

## §1 Preliminaries on the Subgroup Structure of $McL$

The McLaughlin group $McL$, first discovered in 1969 [McL], is a sporadic simple group of order $898,128,000 = 2^7 \cdot 3^6 \cdot 5^3 \cdot 7 \cdot 11$. At the prime $p = 2$ it has rank four i.e. its largest 2–elementary abelian subgroup is $(\mathbb{Z}/2)^4$. It shares its 2–Sylow subgroup with the Mathieu groups $M_{22}$, $M_{23}$ as well as with $\tilde{A}_8$, the unique perfect central extension of $A_8$ by $\mathbb{Z}/2$, and $PSU_4(3) = PSU_4(\mathbb{F}_9)$.

In this section we describe some basic facts about the subgroup structure of $McL$, which will be used in our cohomology calculations. To begin, we describe $H = Syl_2(M_{22})$. We refer to [AM2] for complete details. This group can be described in six distinct but equally useful ways:

- in exactly two ways as a semidirect product $(\mathbb{Z}/2)^4 \times_T D_8$;
- as a central extension $2{\cdot}UT_4(2)$ (where $UT_4(2)$ denotes the group of upper triangular $4 \times 4$ matrices over $\mathbb{F}_2$);
- in exactly two ways as a split extension $UT_4(2){:}\,2$;
- as a split extension $Syl_2(L_3(4)).2$, where the $\mathbb{Z}/2$ action is induced by the Frobenius map.

From now on we write $S$ for the group $S = Syl_2(L_3(4)) \subset H$.

Lemma 1.1:

*(1) There are exactly two $2^4$'s contained in $H$, each normal, so $H = 2^4{:}\,D_8$ in exactly two distinct ways.*

*(2) The two $2^4$'s above are not conjugate in any of the simple groups above and their normalizers are given by the following table*

| | $M_{22}$ | $M_{23}$ | $PSU_4(3)$ | $McL$ |
|---|---|---|---|---|
| $2^4_I$ | $2^4{:}\,\mathcal{A}_6$ | $2^4{:}\,\mathcal{A}_7$ | $2^4{:}\,\mathcal{A}_6$ | $2^4{:}\,\mathcal{A}_7$ |
| $2^4_{II}$ | $2^4{:}\,\mathcal{S}_5$ | $2^4{:}\,GL_2(4){:}\,2$ | $2^4{:}\,\mathcal{A}_6$ | $2^4{:}\,\mathcal{A}_7$ |

*Here the action of the $\mathbb{Z}/2$ on $GL_2(4)$ is via the Galois automorphism of the field $\mathbb{F}_4$.*

*(3) Both $Aut(PSU_4(3))$ and $Aut(McL)$ contain elements of order two which exchange $2^4_I$ and $2^4_{II}$.* ∎

(The first statement is proved in [J]. The second and third statements are contained in [Co], a more detailed description of the normalizer of $2^4_{II}$ for $M_{23}$ is given in [M].)

Remark 1.2: The group $S$ is given as the span $\langle 2^4_I, 2^4_{II}\rangle$. In particular, from 1.1(1) it follows that there is only one copy of $S$ contained in $H$, and any copy of $L_3(4) \subset G$ where $G$ is one of the four simple groups above must intersect some $G$–conjugate of $H$ in $S$. Also, the subgroup $2^4_{II}{:}\,GL_2(4) \subset 2^4_{II}{:}\,GL_2(4){:}\,2 \subset M_{23}$ has $S$ as its 2-Sylow subgroup.

Consider one of the two copies of $UT_4(2) \subset H$, say the one which contains $2^4_I$ for definiteness. It can be written in the form $2^4_I{:}\,W$, where

$$W = \langle \alpha, \beta\rangle \cong 2^2.$$

The action of $W$ on $2^4_I$ satisfies the condition that the intersection of the centralizer for each of the three non-trivial elements in $W$ with $2^4_I$ is a $2^2$, and these three $2^2$'s have a single $\mathbb{Z}/2$ intersection which we write $\langle c\rangle$. Also, $c$ generates the center of $UT_4(2)$ and of $H$. Thus there are three $2^3 \subset UT_4(2)$ each of which intersects $W$ in a single 2 that we call $W_\alpha$, $W_\beta$, and $W_{\alpha\beta}$ for their interesections with $W$.

In $H$, following [AM2], we can assume that $W_\alpha$ and $W_\beta$ are conjugate, but $W_{\alpha\beta} \subset 2^4_{II}$ is not conjugate to the others. Also, there is another $2^3$ subgroup given as $\mathcal{B} = \langle \alpha, \beta, c\rangle$. Under the automorphism of 1.1(3) which exchanges the two $2^4$'s in $H$, the image of $\mathcal{B}$ is $W_\alpha$ or $W_\beta$. In fact we have the following result from [AM2], [M].

Lemma 1.3:

*(1) There are 4 conjugacy classes of extremal 2-elementary subgroups contained in $H$, the two $2^4$'s, and two classes of $2^3$'s, each class containing 2 subgroups, the first $\{W_\alpha, W_\beta\}$ and the second $\{\mathcal{B}, \bar{\mathcal{B}}\}$.*

*(2) In $M_{22}$ there is an element $\gamma \in 2^4_I{:}\,\mathcal{A}_6$ of order three which conjugates $W_\alpha$ to $W_{\alpha\beta}$. Thus, the two extremal $2^3$'s in the first class fuse with a $2^3 \subset 2^4_{II}$ in $M_{22}$.*

*(3) In $M_{22}$ and $M_{23}$ the $2^3$ represented by $\mathcal{B}$ remains extremal.* ∎

We next recall the notion of weak closure.

Definition 1.4: *Given a triple of groups $V \subset K \subset G$, we say that $V$ is weakly closed in $K$ if any subgroup $V' \subset K$ which is conjugate to $V$ in $G$ is already conjugate to $V$ in $K$.* ∎

Now we have

Lemma 1.5: *In $PSU_4(3)$ and $McL$ the two $2^4$'s are weakly closed in $H \subset PSU_4(3)$ and the $2^3$ represented by $\mathcal{B} \subset H$ is conjugate to a subgroup of one of the two $2^4$'s while $W_\alpha$ is conjugate to a subgroup of the other $2^4$. Consequently, the only extremal 2-elementaries in $PSU_4(3)$ are the two $2^4$'s.*

Proof: Since the two $2^4$'s are not conjugate in $PSU_4(3)$ they are weakly closed from the definition. Thus it remains to show that both families of extremal $2^3$'s are conjugate to subgroups of these $2^4$'s.

Consider the normalizers of the two $2^4$'s in $PSU_4(3)$. Each has the form of a split extension $2^4{:}\,\mathcal{A}_6$, hence is completely determined by the action of $\mathcal{A}_6$ on $2^4$, but there is only one conjugacy class of subgroups isomorphic to $\mathcal{A}_6$ contained in $\mathcal{A}_8 = L_4(2)$, so these extensions $2^4{:}\mathcal{A}_6$ are each isomorphic to the group with the same name in $M_{22}$. Consequently, the first identifies $W_\alpha$ with a subgroup of the other $2^4$ exactly as in $M_{22}$.

One of the outer automorphisms of $PSU_4(3)$ when restricted to $H$ gives the automorphism above exchanging the two $2^4$'s and the families of extremal $2^3$. It follows that $\mathcal{B}$ must be conjugate to a subgroup of the other $2^4$ in $PSU_4(3)$ as asserted. ∎

There is a sporadic geometry discovered by Ronan associated to the group $McL$ (see [RSY] for details). From the topological point of view, this means that there is a finite CW–complex $X$ with an action of $McL$ arising from its internal combinatorial structure. Further, in [RSY] it was shown that for this complex the fixed point set under any element of order two in $McL$ is contractible, and hence its stabilizers can be used to compute $H^*(McL)$. However, as this would involve extensive computations we have chosen instead to make a calculation using more direct methods. Nevertheless it contains essential information about the subgroup structure of $McL$. The orbit space is a single solid 3–dimensional simplex, described by the following table:

| Simplex | Stabilizer | Simplex | Stabilizer |
|---|---|---|---|
| $a$ | $PSU_4(3)$ | $bc$ | $2^4 : (3 \times A_4) : 2$ |
| $b, c$ | $2^4 : A_7$ | $bd, cd$ | $2^4 : L_3(2)$ |
| $d$ | $\tilde{A}_8$ | $abd, acd, bcd$ | $2^4 : S_4$ |
| $ab, ac$ | $2^4 : A_6$ | $abc$ | $2^4 : S_4$ |
| $ad$ | $2^{1+4} : (S_3 \times S_3)$ | $abcd$ | $2^4 : D_8$ |

We can best visualize this geometry by considering the following figure, where the larger type identifies the isotropy groups of the vertices, the smaller type gives the isotropy groups of the edges, each face has isotropy group $2^4{:}\mathcal{S}_4$, and the interior has isotropy group $H$.

1.6

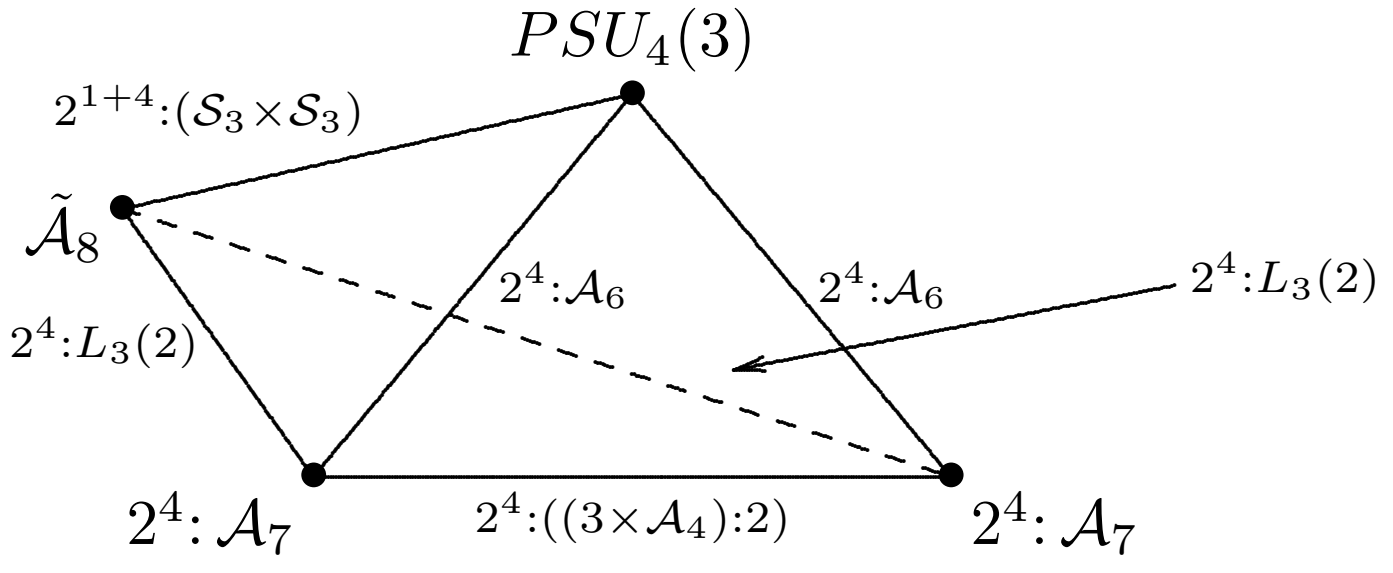

THE ORBIT SPACE OF THE SPORADIC GEOMETRY FOR $McL$

Actually in [RSY] it is also shown that there is a sub–geometry corresponding to the bottom face of the simplex above from which the cohomology of $McL$ can be determined.

This 2-dimensional geometry is similar to the one associated to $M_{22}$ (see [AM2]). In fact, that group is known to be isomorphic to the generalized amalgam of the stabilizers. Whether or not an analogue of this result is true for $McL$ remains to be verified.

From the above, the group $2^4{:}\,(3 \times \mathcal{A}_4).2$ is the intersection of the two $2^4{:}\,\mathcal{A}_7$'s which normalize $2^4_I$ and $2^4_{II}$ respectively. Since $McL$ contains $M_{22}$ we see that $Syl_2(2^4 : 3 \times \mathcal{A}_4) = S$, and that there are two elements of order three, denoted $T_1$ and $T_2$ acting on $S$ which gives an embedding $3^2 \subset \mathrm{Out}(S)$.

Lemma 1.7: *There is an inclusion* $S \subset 2^4_{II}{:}\,\mathcal{S}_5 \subset M_{22}$ *as the inverse image of the Klein group* $K \subset \mathcal{S}_4 \subset \mathcal{S}_5$.

Proof: We know that in $M_{22}$ we have

$$2^4_{II}{:}\,\mathcal{S}_5 \cap 2^4_I{:}\,\mathcal{A}_6 \;=\; 2^4_{II}{:}\,\mathcal{S}_4$$

where the first two groups are the normalizers of the two distinct $2^4$'s in $Syl_2(M_{22}) = 2^{2+4}{:}\,2$. There is only one conjugacy class of subgroups $\mathcal{S}_4$ in $\mathcal{S}_5$ and the intersection of the second $2^4$ with it must be a normal subgroup, hence $K$ as asserted. ∎

Corollary 1.8: *$Aut(S)$ has the form*

$$V{:}\,\mathcal{S}_3 \times \mathcal{S}_3$$

*where $V$ is a finite 2-group.*

Proof: Consider the sequence of inclusions

$$2^4{:}\,\mathcal{S}_4 \hookrightarrow 2^4{:}\,\mathcal{S}_5 \hookrightarrow 2^4{:}\,\mathcal{A}_7 \hookrightarrow 2^4{:}\,\mathcal{A}_8$$

where $\mathcal{A}_8 = L_4(2)$ and, as in 1.7, $S = 2^4{:}\,K \subset 2^4{:}\,\mathcal{S}_4$. It is direct to check that if an automorphism, $\lambda$, of $S$ is the identity on $2^4_{II}$ then it is given by specifying a homomorphism $\phi{:}\,K \to C \subset 2^4_{II}$ where $C = 2^2 \subset 2^4_{II}$ is the center of $S$, and setting $\lambda(k) = k + \phi(k)$ for each $k \in K$. Clearly, this subgroup $V \subset \mathrm{Aut}(S)$ is normal in the subgroup $M \subset \mathrm{Aut}(S)$ which normalizes $2^4_{II}$, and, since $V = 2^4$, we have

$$M \;=\; 2^4 \cdot N_{\mathcal{A}_8}(K)$$

where $K \subset \mathcal{S}_4$ is the Klein group. Also, the entirety of the Aut-group is the extension of $M$ by the 2 which switches the two $2^4$'s. But we have easily that

$$N_{\mathcal{A}_8}(K) = \mathcal{A}_4 \times \mathcal{A}_4{:}\,2$$

and from this the result is a direct calculation. ∎

## §2 Preliminaries on the cohomology of $PSL_4(3)$ and $McL$

From the Quillen-Venkov detection theorem [AM1, p.144] we have an exact sequence for any finite group $G$,

$$0 \longrightarrow \mathrm{rad}(H^*(G)) \longrightarrow H^*(G) \xrightarrow{res^*} \coprod_V H^*(V)^{W_G(V)},$$

where the $V$ run over representatives for the conjugacy classes of maximal 2-elementary subgroups of $G$ and $W_G(V) = N_G(V)/V$ is the Weyl group of $V$ in $G$.

In [AM2] we proved that the Lyndon–Hochschild–Serre spectral sequence associated to writing $H$ as the extension $S.2$ collapses at the $E_2$–stage, yielding a precise description of $H^*(H)$. It is most convenient to describe the result there in the following manner. First, the radical is given as

2.1 $$\mathrm{rad}(H^*(H)) \cong \mathrm{rad}(H^*(S))^{\mathbb{Z}/2} = \mathbb{F}_2[t_4, t_8](a_2, \alpha_3, T_6, T_7).$$

Here, for the sake of clarity, we emphasize that the radical of $H^*(H)$ injects to $\mathrm{rad}(H^*(S))$ under the restriction map. Next, the image of $H^*(H)$ in each of the four extremal $2^n$'s given in 1.2 is described. In particular, *the image in each* $2^4$ *is the full invariant subalgebra* $H^*(2^4)^{D_8}$. (The image in each of the two conjugacy classes of extremal $2^3$'s is a proper subset of its invariant algebra, however in the current situation, as they fuse with subgroups of the two $2^4$'s, this need not concern us.)

LEMMA 2.2:
*(a)* *rad*$(H^*(PSU_4(3))) \subset$ *rad*$(H^*(M_{22})) \cong \mathbb{F}_2[d_8, d_{12}](a_2, a_7, a_{11}, a_{14})$.
*(b)* *rad*$(H^*(McL)) \subset$ *rad*$(H^*(M_{23})) = \mathbb{F}_2[d_8, d_{12}](a_7, a_{11})$.

PROOF: From [Co] we have that $PSU_4(3) \subset McL$ contains a subgroup $L_3(4)$. Consequently, it follows from 2.1 and 1.2 that

$$\mathrm{rad}(H^*(PSU_4(3))) \subset \mathrm{rad}(H^*(L_3(4))) \cap H^*(S)^{\mathbb{Z}/2}.$$

But from [AM2], the intersection above is $\mathrm{rad}(H^*(M_{22})) = \mathbb{F}_2[d_8, d_{12}](a_2, a_7, a_{11}, a_{14})$.

In [M] the radical of $H^*(M_{23})$ was shown to be

$$\mathrm{rad}(M_{22})^{\mathbb{Z}/3} = \mathbb{F}_2[d_8, d_{12}](a_7, a_{11})$$

where, from 1.1, the 3-Sylow subgroup of $N_{M_{23}}(S)$ is equal to $(\mathbb{Z}/3)^2$ with the first $\mathbb{Z}/3 \subset L_3(4)$ and the second giving the action above. Consequently, from 1.7, 1.8, we see that $\mathrm{rad}(H^*(McL))$ is also invariant under the same $\mathbb{Z}/3$-action and 2.2 follows. ∎

It remains to describe the exact structure of these radicals, and to describe the restriction of $H^*(PSU_4(3))$ and $H^*(McL)$ to the cohomology of the two $2^4$'s. It turns out that the function of the radicals is to provide elements with higher Bocksteins to make sure that all the elements in $H^*(G)$ are torsion where $G$ is $PSU_4(3)$ or $McL$. Consequently, we will be able to show that a sufficient number of elements are in the radicals in each case once we have some information on the image of restriction in $H^*(2^4)$.

Recall the Cardenas-Kuhn theorem [AM1, p.113]:

THEOREM 2.3: *Let $L \subsetneq K \subsetneq G$ be weakly closed with $L$ a $p$-elementary subgroup of $G$. Suppose that $N_K(L)$ contains a $p$-Sylow subgroup of $N_G(L)$, then the image of*

$$res^*: H^*(G; \mathbb{F}_p) \longrightarrow H^*(L; \mathbb{F}_p)$$

*is exactly equal to the intersection*

$$H^*(L; \mathbb{F}_p)^{W_G(L)} \cap im(res^*: H^*(K; \mathbb{F}_p) \longrightarrow H^*(L; \mathbb{F}_p)).$$

∎

Thus we have, from 1.5 and the remarks before 2.1

COROLLARY 2.4: *The restriction image from $H^*(PSU_4(3)) \subset H^*(2^4_J)$ is the full invariant algebra $H^*(2^4)^{\mathcal{A}_6}$ for $J = I, II$. Similarly the restriction image from $H^*(McL) \subset H^*(2^4_J)$ is the full invariant algebra $H^*(2^4)^{\mathcal{A}_7}$ for $J = I, II$.* ∎

From [AM4] we have that

2.5 $$H^*(2^4)^{\mathcal{A}_6} = \mathbb{F}_2[\bar{w}_3, \gamma_5, d_8, d_{12}](1, \gamma_9, b_{15}, \gamma_9 b_{15})$$

and

2.6 $$H^*(2^4)^{\mathcal{A}_7} = \mathbb{F}_2[d_8, d_{12}, d_{14}, d_{15}](1, x_{18}, x_{20}, x_{21}, x_{24}, x_{25}, x_{27}, x_{45})$$

which specifies these invariant subalgebras. Consequently, in order to specify the rings $H^*(PSU_4(3))$ and $H^*(McL)$ we need to finish the determination of the radicals in each case and the image in

$$H^*(2^4_I)^{W_G(2^4)} \oplus H^*(2^4_{II})^{W_G(2^4)}$$

which really means specifying which elements $(a, 0)$ are in the restriction image for $a$ in the invariant subalgebra. We know from our analysis of $H^*(H)$ in [AM2] that the element $(\bar{w}_3, 0)$ is not in the image of restriction and neither are $(d_8, 0)$, $(d_{12}, 0)$. In fact, the *quotient algebra* of elements not of this form must contain at least

$$\mathbb{F}_2[d_8, d_{12}](1, \bar{w}_3).$$

In any case, $a_2 \in$ rad and $(\bar{w}_3, \bar{w}_3)$ are certainly elements in $H^*(PSU_4(3))$. In fact it will turn out that $\text{rad}(H^*(PSU_4(3))) \cong \text{rad}(H^*(L_3(4)) \cap H^*(S)^{\mathbb{Z}/2}$ by a fairly easy argument. But the determination of which elements in $H^*(E_1)^{\mathcal{A}_6} \oplus H^*(E_2)^{\mathcal{A}_6}$ are not doubled in the image from $H^*(PSU_4(3))$ will necessitate an independent calculation of $H^*(PSU_4(3))$ which we turn to next.

## §3 The Structure of $H^*(PSU_4(3))$

The procedure we use in this section to determine the ring $H^*(PSU_4(3))$ is to use our knowledge of the ring $H^*(SU_4(3))$ (see [FP]) and an Eilenberg-Moore spectral sequence with $E_2$-term $\mathrm{Ext}_{H_*(C(SU_4(3)))}(H_*(SU_4(3)), \mathbb{F}_2)$ that converges to $H^*(PSU_4(3))$ to make initial calculations (here $C(SU_4(3))$ denotes the center of $SU_4(3)$). Then, coupling these initial results with the discussion in the previous section will enable us to determine $H^*(PSU_4(3))$ completely.

*The Eilenberg-Moore spectral sequence for* $PSU_4(3)$

We know $H^*(SU_4(3)) = \mathbb{F}_2[b_4, b_6, b_8] \otimes E(e_3, e_5, e_7)$ from [FP]. We also have that $SU_4(3)$ is a central extension of $PSU_4(3)$, $C \triangleleft SU_4(3) \to PSU_4(3)$ with center $C = \mathbb{Z}/4 \subset \mathbb{F}_9^*$. Each central extension $C \overset{\triangleleft}{\to} G \to G/C$ is classified by an element $\alpha \in H^2(B_{G/C}; C)$, which is the same thing as a map $\alpha\colon B_{G/C} \to K(C, 2)$ with $K(C, 2)$ an Eilenberg-MacLane space. That is to say $K(C, 2)$ is a CW-complex with homotopy groups given as follows:

$$\pi_i(K(C,2)) = \begin{cases} C & \text{if } i = 2, \\ 0 & \text{if } i \neq 2, \end{cases}$$

and it is well known that $H^2(X; C)$ is naturally identified with the set of homotopy classes of maps $[X, K(C, 2)]$ for $X$ any CW-complex. In fact, in our situation, if the map $\alpha\colon B_{PSU_4(3)} \to K(C, 2)$ is chosen to represent the $K$-invariant of the extension, then there is a fibration sequence for the classifying spaces

$$B_{\mathbb{Z}/4} \longrightarrow B_{SU_4(3)} \longrightarrow B_{PSU_4(3)} \overset{\alpha}{\longrightarrow} K(\mathbb{Z}/4, 2)$$

and the leftmost fibration is principal.

There are many spectral sequences which can be used to determine the cohomology of the total space of a fibration $F \to E \to B$. The Serre spectral sequence is the most commonly used, with $E_2$-term $H^*(B; H^*(F))$. However, when we know the fiber and total space, but want the cohomology of the base then the most efficient tool is the Eilenberg-Moore spectral sequence provided the fibration is principal. Its $E_2$-term is equal to

$$\mathrm{Ext}_{H_*(F;\mathbb{F})}(H_*(E; \mathbb{F}), \mathbb{F})$$

and it converges to $H^*(B; \mathbb{F})$ for $\mathbb{F}$ any field.

To determine this $E_2$-term we need to study the action of $H_*(B_{\mathbb{Z}/4})$ on $H_*(B_{SU_4(3)})$. In fact this action arises from the *homomorphism* $\mathbb{Z}/4 \times SU_4(3) \to SU_4(3)$ obtained by multiplication by the center, on passing to classifying spaces,

$$(*) \qquad B(\times)\colon B_{\mathbb{Z}/4} \times B_{SU_4(3)} = B_{\mathbb{Z}/4 \times SU_4(3)} \longrightarrow B_{SU_4(3)}.$$

Write $H^*(\mathbb{Z}/4) = \mathbb{F}_2[d] \otimes E(e)$ where $e$ is the 1-dimensional generator, and $d$ is the non-zero element in dimension 2. Then we have

Lemma 3.1: *The cohomology map $H^*(SU_4(3); \mathbb{F}_2) \longrightarrow H^*(\mathbb{Z}/4; \mathbb{F}_2) \otimes H^*(SU_4(3); \mathbb{F}_2)$ associated to $(*)$ above is given on generators by $b_i \mapsto 1 \otimes b_i$, $x_i \mapsto 1 \otimes x_i$ for $i = 3, 4, 5, 6$ but*

$$\begin{aligned} x_7 &\mapsto x_7 + d^2 \otimes x_3 + d \otimes x_5 + e \otimes b_6, \\ b_8 &\mapsto d^4 \otimes 1 + d^2 \otimes b_4 + d \otimes b_6 + 1 \otimes b_8. \end{aligned}$$

Proof: It is easiest to work with $U_4(3)$ with maximal torus $(\mathbb{Z}/4)^4$. The restriction of $H^*(U_4(3); \mathbb{F}_2) = \mathbb{F}_2[b_2, b_4, b_6, b_8] \otimes E(x_1, x_3, x_5, x_7)$ is injective to

$$\mathbb{F}_2[\sigma_1, \sigma_2, \sigma_3, \sigma_4] \otimes E$$

where $E$ is an exterior algebra with generators

$$\begin{aligned} &e_1 + e_2 + e_3 + e_4, \\ &e_1(d_2 + d_3 + d_4) + e_2(d_1 + d_3 + d_4) + e_3(d_1 + d_2 + d_4) + e_4(d_1 + d_2 + d_3), \\ &e_1(\sigma_2(d_2, d_3, d_4)) + e_2(\sigma_2(d_1, d_3, d_4) + e_3(\sigma_2(d_1, d_2, d_4)) + e_4(\sigma_2(d_1, d_2, d_3)) \\ &e_1 d_2 d_3 d_4 + e_2 d_1 d_3 d_4 + e_3 d_1 d_2 d_4 + e_4 d_1 d_2 d_3 \end{aligned}$$

and $\sigma_i$ is the $i^{th}$ symmetric function in the $d_i$'s ([FP], pg. 288).

The inclusion $SU_4(3) \subset U_4(3)$ induces a surjection in cohomology onto

$$H^*(SU_4(3); \mathbb{F}_2) = \mathbb{F}_2[b_4, b_6, b_8] \otimes E(x_3, x_5, x_7)$$

where $b_{2i}$ is the image under restriction of $b_{2i}$, and $x_{2i+1}$ is the image under restriction of $x_{2i+1}$. More exactly, the inclusion of the maximal torus for $SU_4(3)$ can be given as $t_i \mapsto t_i t_4^{-1}$, $1 \leq i \leq 3$, so in cohomology $d_i \mapsto d_i$, $e_i \mapsto e_i$, $1 \leq i \leq 3$, while $e_4 \mapsto e_1 + e_2 + e_3$, $d_4 \mapsto d_1 + d_2 + d_3$. Consequently, the images of the generators for $H^*(SU_4(3); \mathbb{F}_2)$ in the cohomology of the maximal torus are given as

$$\begin{aligned} b_4 &\mapsto \sigma_2 + \sigma_1^2 \\ &\quad = d_1^2 + d_2^2 + d_3^2 + d_1 d_2 + d_1 d_3 + d_2 d_3 \\ b_6 &\mapsto \sigma_3 + \sigma_1 \sigma_2 \\ &\quad = d_1 d_2(d_1 + d_2) + d_1 d_3(d_1 + d_3) + d_2 d_3(d_2 + d_3) \\ b_8 &\mapsto \sigma_3 \sigma_1 \\ x_3 &\mapsto e_1(d_2 + d_3) + e_2(d_1 + d_3) + e_3(d_1 + d_2) \\ x_5 &\mapsto e_1(d_2^2 + d_3^2) + e_2(d_1^2 + d_3^2) + e_3(d_1^2 + d_2^2) \\ x_7 &\mapsto e_1 d_2 d_3(d_2 + d_3) + e_2 d_1 d_3(d_1 + d_3) + e_3 d_1 d_2(d_1 + d_2). \end{aligned}$$

We use this to figure out the cohomology map $H^*(SU_3(4)) \longrightarrow H^*(\mathbb{Z}/4) \otimes H^*(SU_3(4))$. Restricted to the maximal torus, the central multiplication

$$\mathbb{Z}/4 \times (\mathbb{Z}/4)^3 \longrightarrow (\mathbb{Z}/4)^3$$

is given by $(\tau, t_1, \ldots, t_3) \mapsto (\tau t_1, \tau t_2, \tau t_3)$ and consequently, in cohomology we have

$$\begin{aligned} e_i &\mapsto e_i + e, \ 1 \leq i \leq 3, \\ d_i &\mapsto d_i + d, \ 1 \leq i \leq 3, \end{aligned}$$

so $b_4 \mapsto 1 \otimes b_4$, $b_6 \mapsto 1 \otimes b_6$, but

$$\begin{aligned} b_8 &\mapsto (\sigma_1 + d)(\sigma_3 + d\sigma_2 + d^2\sigma_1 + d^3) \\ &= \sigma_1\sigma_3 + d(\sigma_3 + \sigma_1\sigma_2) + d^2(\sigma_1^2 + \sigma_2) + d^4 \\ &= b_8 + d \otimes b_6 + d^2 \otimes b_4 + d^4 \otimes 1. \end{aligned}$$

The calculation for the odd dimensional generators is similar. ∎

As a consequence

$$\Delta b_8^{2^n} = d^{2^{n+2}} + d^{2^{n+1}} b_4^{2^n} + d^{2^n} b_6^{2^n} + b_8^{2^n}$$

and $H^*(SU_4(3))$ is free over the quotient coalgebra $\mathbb{F}_2[d^4]$. (Dually, $H_*(SU_4(3))$ is free over the subalgebra $E((d^4)^*, (d^8)^*, \cdots)$ where $H_*(\mathbb{Z}/4; \mathbb{F}_2) = E(e, d^*, (d^2)^*, (d^4)^*, \ldots)$, with generators the dual elements to the cohomology subalgebra $\mathbb{F}_2[b_4, b_6] \otimes E(x_3, x_5, x_7)$.)

Now that we have the module structure of $H_*(SU_4(3))$ over $H_*(\mathbb{Z}/4)$ (or dually, the comodule structure of $H^*(SU_4(3))$ over $H^*(\mathbb{Z}/4)$) we are ready to determine the groups $\mathrm{Ext}_{H_*(\mathbb{Z}/4)}(H_*(SU_4(3)), \mathbb{F}_2)$. The procedure will be to use the freeness over the subalgebra above and the change of rings theorem:

THEOREM 3.2 (Change of rings): *Let $\mathcal{A}$ be an augmented Hopf algebra over a field $\mathbb{F}$ and $\mathcal{M}$ an $\mathcal{A}$-module which is free over the sub-Hopf-algebra $\mathcal{B} \subset \mathcal{A}$. Let $\mathcal{C} = \mathcal{A}//\mathcal{B}$ be the quotient Hopf algebra and $\mathcal{M}//\mathcal{B}$ the quotient module, then*

$$Ext_{\mathcal{A}}(\mathcal{M}, \mathbb{F}) \cong Ext_{\mathcal{C}}(\mathcal{M}//\mathcal{B}, \mathbb{F}).$$

∎

In the present case let us denote the quotient Hopf algebra

$$H_*(\mathbb{Z}/4)//E((d^4)^*, (d^8)^*, \ldots) = E(e, d^*, (d^2)^*) \text{ as } \Lambda,$$

and the quotient module $H_*(SU_4(3))/\mathcal{I}(E((d^4)^*, (d^8)^*, \ldots)$ as $\mathcal{N}$. In $\mathcal{N}$ the action of $\Lambda$ is given by $e(b_6^j)^* = (b_6^{j-1}x_7)^*$, $d(x_5\theta)^* = (x_7\theta)^*$, $d^2(x_3\tau)^* = (x_7\tau)^*$ and otherwise the action is trivial.

Using change of rings we have

$$\begin{aligned}\mathrm{Ext}_{H_*(\mathbb{Z}/4)}(H_*(SU_4(3)),\mathbb{F}_2) &\cong \mathrm{Ext}_\Lambda(\mathcal{N},\mathbb{F}_2)\\ &\cong \mathrm{Ext}_\Lambda(\mathbb{F}_2[b_4,b_6]\otimes E(x_3,x_5,x_7)^*,\mathbb{F}_2)\\ &\cong \mathbb{F}_2[b_4]\otimes \mathrm{Ext}_\Lambda(\mathbb{F}_2[b_6]\otimes E(x_3,x_5,x_7)^*,\mathbb{F}_2),\end{aligned}$$

since the action of $\Lambda$ does not involve $\mathbb{F}_2[b_4]$. Additionally, since

$$\mathrm{Ext}_{A\otimes B}(\mathbb{F},\mathbb{F}) = \mathrm{Ext}_A(\mathbb{F},\mathbb{F})\otimes \mathrm{Ext}_B(\mathbb{F},\mathbb{F}),$$

and $\mathrm{Ext}_{E(y)}(\mathbb{F}_2,\mathbb{F}_2) = \mathbb{F}_2[|y|]^\dagger$ we have that $\mathrm{Ext}_\Lambda(\mathbb{F}_2,\mathbb{F}_2) = \mathbb{F}_2[\lambda_3,\lambda_5,\lambda_7]$, the polynomial ring on the stated generators where $\lambda_i \in \mathrm{Ext}_\Lambda^{i-1,1}(\mathbb{F}_2,\mathbb{F}_2)$.

Before we determine the Ext-groups $\mathrm{Ext}_\Lambda(\mathbb{F}_2[b_4,b_6]\otimes E(x_3,x_5,x_7)^*,\mathbb{F}_2)$ we point out that they form an associative, commutative, ring with unit – indeed, a ring over $\mathrm{Ext}_\Lambda(\mathbb{F}_2,\mathbb{F}_2)$. To see this, note that the diagonal map gives a ring homomorphism

$$H_*(\mathbb{Z}/4)\longrightarrow H_*((\mathbb{Z}/4)^2)\cong H_*(\mathbb{Z}/4)\otimes H_*(\mathbb{Z}/4)$$

and the diagonal map for $SU_4(3)$ gives a module map

$$H_*(SU_4(3))\longrightarrow H_*(SU_4(3))\otimes H_*(SU_4(3)).$$

Together these two maps thus give a pairing map

$$\mathrm{Ext}_{H_*(\mathbb{Z}/4)}(H_*(SU_4(3)),\mathbb{F}_2)\otimes \mathrm{Ext}_{H_*(\mathbb{Z}/4)}(H_*(SU_4(3)),\mathbb{F}_2)\longrightarrow \mathrm{Ext}_{H_*(\mathbb{Z}/4)}(H_*(SU_4(3),\mathbb{F}_2)$$

which has the requisite properties. In particular, the Ext groups have the form of the tensor product

$$\mathbb{F}_2[b_4]\otimes \mathrm{Ext}_\Lambda(\mathbb{F}_2[b_6]\otimes E(x_3,x_5,x_7)^*,\mathbb{F}_2)$$

as a ring.

Theorem 3.3: *We have*

$$Ext_\Lambda(\mathbb{F}_2[b_6]\otimes E(x_3,x_5,x_7)^*,\mathbb{F}_2)\cong \mathbb{F}_2[b_6,x_3,x_5,\lambda_2,\lambda_3,\lambda_5]/\mathcal{R}$$

*where $\mathcal{R}$ is the ideal with generators*

$$\mathcal{R} = (x_3^2, x_5^2, \lambda_2 b_6 + \lambda_3 x_5 + \lambda_5 x_3).$$

† Here $E(y)$ is a regarded as a primitively generated Hopf algebra so that $\mathrm{Ext}_{E(y)}(\mathbb{F},\mathbb{F})$ becomes a commutative ring.

REMARK 3.4: The following figure gives the groups $\mathrm{Ext}_\Lambda(\mathbb{F}_2[b_6] \otimes E(x_3, x_5, x_7)^*, \mathbb{F}_2)$ through dimension 8, the range which is most important to calculate $H^*(PSU_4(3))$.

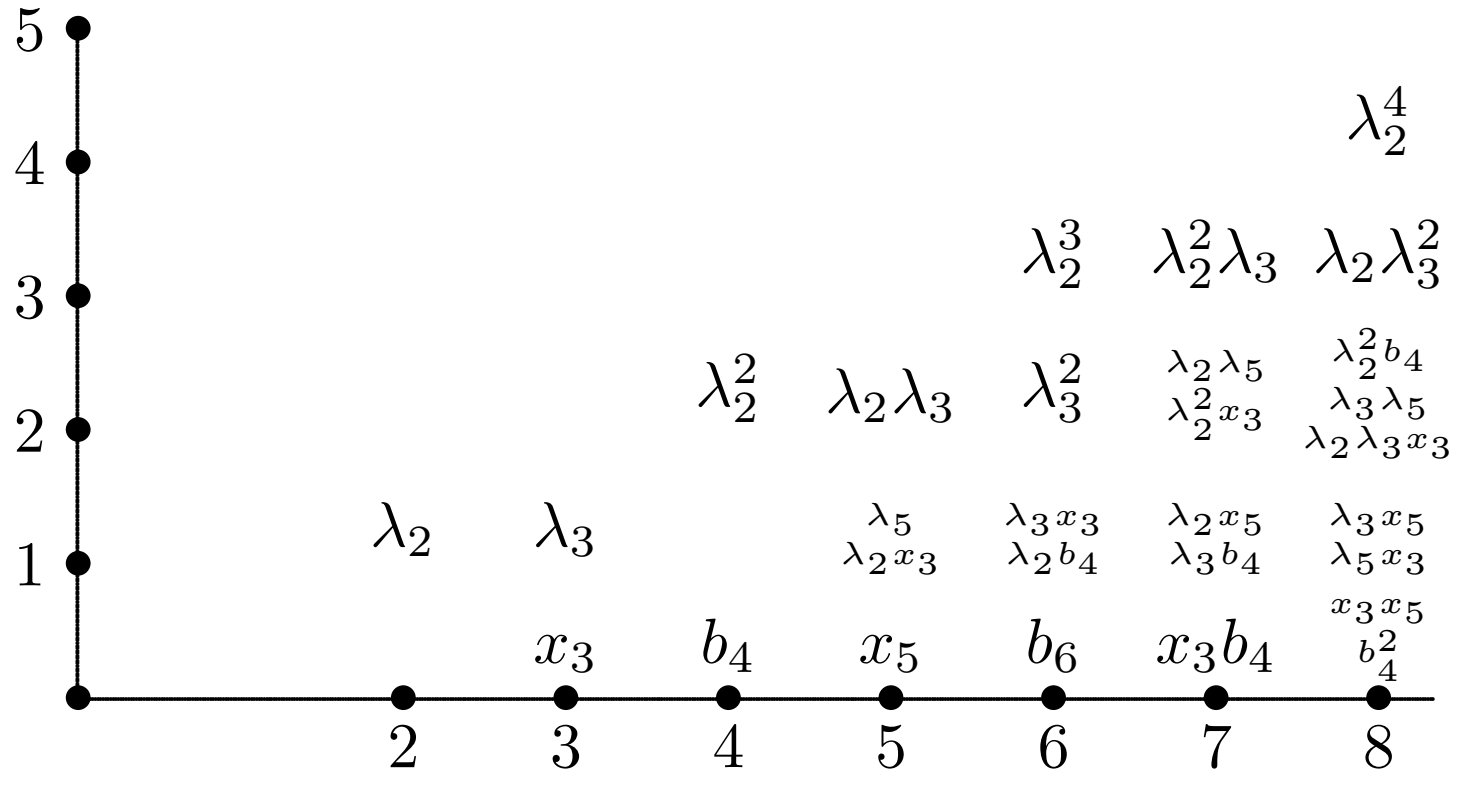

**The $E_2$-term of the Eilenberg-Moore spectral sequence for $PSU_4(3)$**

PROOF: Since $e(b_6^j)^* = (b_6^{j-1}x_7)^*$ for $j \geq 1$ it follows that

$$\mathcal{F}_i = E(x_3, x_5, x_7)(1, b_6, b_6^2, b_6^3, \ldots, b_6^i)^*$$

is a sub-$\Lambda$ module of $\mathbb{F}_2[b_6] \otimes E(x_3, x_5, x_7)^*$ for each $i = 1, 2, \ldots$. This defines a filtration of $\mathbb{F}_2[b_6] \otimes E(x_3, x_5, x_7)^*$, and the resulting spectral sequence has $E_1^i$-term equal to $\mathrm{Ext}_\Lambda(E(x_3, x_5, x_7)^*, \mathbb{F}_2)b_6^i$.

We now determine this Ext-group.

LEMMA 3.5:

$$Ext_\Lambda(E(x_3, x_5, x_7)^*, \mathbb{F}_2) \cong \mathbb{F}_2[\lambda_2, \lambda_3, \lambda_5](1, x_3, x_5, x_3x_5, x_{15}, y_{15})/\mathcal{R}$$

*where* $y_{15} \in Ext_\Lambda^{14,1}$, $x_{15} \in Ext^{15,0}$ *and* $\mathcal{R}$ *is the ideal*

$$(x_3^2, x_5^2, x_{15}^2, y_{15}^2, \lambda_5x_3 + \lambda_3x_5).$$

*The multiplication is given by the relations in* $\mathcal{R}$ *above and* $x_3y_{15} = \lambda_3x_{15}$, $x_5y_{15} = \lambda_5x_{15}$.

PROOF: We write out the module as follows

$$1, x_3, x_5, x_7, x_3x_5, x_3x_7, x_5x_7, x_3x_5x_7,$$

where $d^*x_5 = (d^2)^*x_3 = x_7$ while $d^*x_3x_5 = x_3x_7$, $(d^2)^*x_3x_5 = x_5x_7$ give the complete action. It follows that $\mathcal{V} = \{x_7, x_3x_7, x_5x_7, x_3x_5x_7\}$ is a trivial $\Lambda$-submodule with quotient $\mathcal{W} = \{x_3, x_5, x_3x_5\}$ also a trivial module. From this short exact sequence of modules we get a long exact sequence of Ext-groups:

$$\cdots \longrightarrow \mathrm{Ext}_\Lambda^*(\mathcal{W}, \mathbb{F}_2) \longrightarrow \mathrm{Ext}_\Lambda^*(E(x_3, x_5, x_7)^*, \mathbb{F}_2) \longrightarrow \mathrm{Ext}_\Lambda^*(\mathcal{V}, \mathbb{F}_2) \xrightarrow{\delta} \mathrm{Ext}_\Lambda^{*+1}(\mathcal{W}, \mathbb{F}_2) \cdots$$

and $\delta(x_7) = \lambda_3 x_5 + \lambda_5 x_3$, $\delta(x_3 x_7) = \lambda_3 x_3 x_5$, $\delta(x_5 x_7) = \lambda_5 x_3 x_5$, while $\delta(\theta x) = \theta\delta(x)$ for $\theta \in \mathrm{Ext}_\Lambda(\mathbb{F}_2, \mathbb{F}_2) = \mathbb{F}_2[\lambda_2, \lambda_3, \lambda_5]$. Also, note that

$$\begin{aligned}\mathrm{Ext}_\Lambda(\mathcal{W}, \mathbb{F}_2) &\cong \mathbb{F}_2[\lambda_2, \lambda_3, \lambda_5] \otimes \mathcal{W} \\ \mathrm{Ext}_\Lambda(\mathcal{V}, \mathbb{F}_2) &\cong \mathbb{F}_2[\lambda_2, \lambda_3, \lambda_5] \otimes \mathcal{V},\end{aligned}$$

as $\mathrm{Ext}_\Lambda(\mathbb{F}_2, \mathbb{F}_2)$-modules, since both $\mathcal{V}$ and $\mathcal{W}$ are trivial over $\Lambda$. From this the calculation follows directly. ∎

Now, to complete the proof of the theorem we need to determine the differentials in the spectral sequence. Since $(d)^*(b_6^n)^* = (x_7 b_6^{n-1})^*$, as has been already pointed out, we see that $d_1(x_{15} b_6^{i-1}) = \lambda_2 x_3 x_5 b_6^i$, and similarly, $d_1(y_{15} b_6^{i-2}) = \lambda_2^2 b_6^i$. Note that both of these differentials are injective on $\mathbb{F}_2[\lambda_2, \lambda_3, \lambda_5]$ tensored with the respective generators.

From this the $E_2$-term is the algebra specified, and it is direct to check that each element is an infinite cycle. ∎

This completes the first stage of the calculation of $H^*(PSU_4(3))$.

*The $d_2$-differentials in the Eilenberg-Moore spectral sequence*

Now that we have determined the $E_2$ term of the Eilenberg-Moore spectral sequence converging to $H^*(PSU_4(3))$ we need to calculate the $d^2$-differentials. We have

Lemma 3.6: *The $d^2$ differentials in the Eilenberg-Moore spectral sequence converging to $H^*(PSU_4(3))$ are completely determined by $d^2(x_3) = \lambda_2^2$, $d^2(b_4) = \lambda_2\lambda_3$, $d^2(x_5) = 0$, $d^2(b_6) = \lambda_2\lambda_5$.*

Proof: Since the spectral sequence is multiplicative and these are the generators it remains only to check that the stated differentials are correct. The first differential is forced since we know from our results on $H^*(M_{22})$ that $a_2^2 = 0$ and $\lambda_2$ must represent $a_2$ in the spectral sequence. Since $w_3 \neq 0$ it must be represented by $\lambda_3$ which appears in filtration $\mathrm{Ext}^{2,1}$ and so cannot have a $d^2$ differential to $\mathrm{Ext}^{2,2}$ (it's image under $d^2$ lies in $\mathrm{Ext}^{1,3} = 0$). Thus, the only element available to truncate $\lambda_2$ is $x_3$. Also, this forces $d^2(b_6) = \lambda_2\lambda_5$, since $x_5$, if it had a $d^2$ differential, could only hit $\lambda_3^2$, a class which is non-zero in $H^*(2^4)^{\mathcal{A}_6}$, and consequently, in $H^*(PSU_4(3))$. It follows, since we have the relation $\lambda_5 x_3 + \lambda_3 x_5 + \lambda_2 b_6 = 0$, that, on passing to differentials $\lambda_5\lambda_2^2 + \lambda_2 d^2(b_6) = 0$, which, as asserted, forces $d^2(b_6) = \lambda_2\lambda_5$.

Now we consider the class $\lambda_2\lambda_3$. Since $\lambda_2$ represents $a_2 \in \mathrm{rad}(H^*(Syl_2(PSU_4(3))))$, and since $L_3(4) \subset PSU_4(3)$ from [Co], it follows that

$$\mathrm{rad}(H^*(PSU_4(3))) \subset \mathrm{rad}(H^*(L_3(4))) \cap H^*(S)^{\mathbb{Z}/2} = \mathbb{F}_2[d_8, d_{12}](a_2, a_7, a_{11}, a_{14}),$$

so $\lambda_3\lambda_2 = 0$ in $H^*(PSU_4(3))$. But the only class which can hit this class is $b_4$. ∎

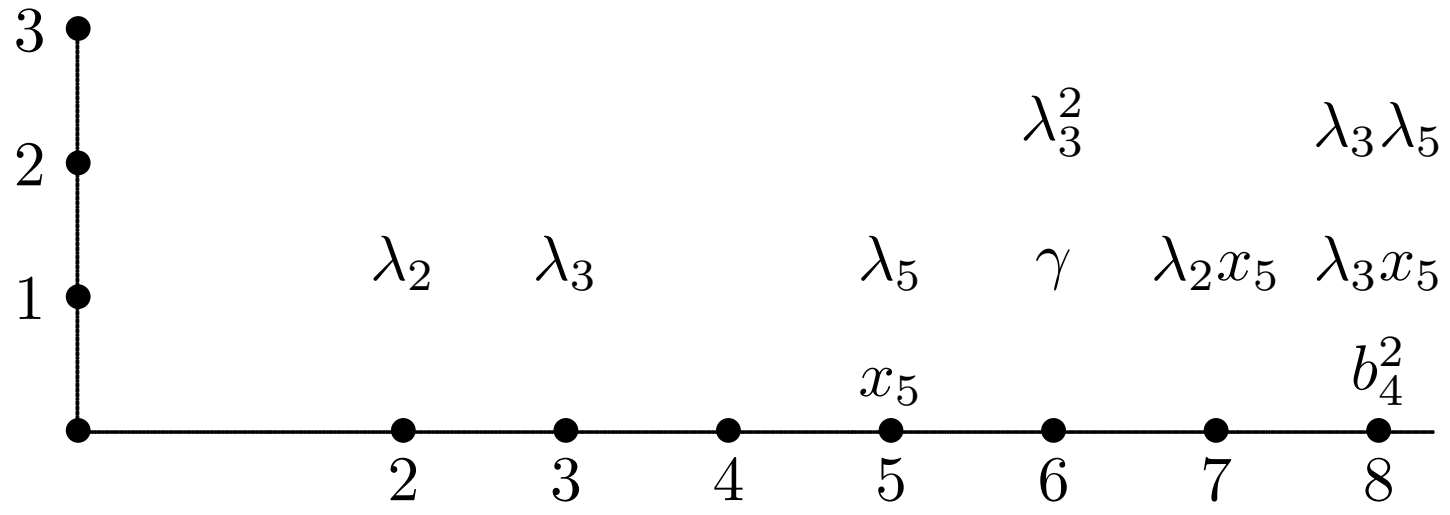

**The $E_3$-term of the Eilenberg-Moore spectral sequence for $PSU_4(3)$**

NOTES: $\lambda_3$ represents $\bar{w}_3$ and $\lambda_5$ represents $Sq^2(\bar{w}_3)$ so $x_5$ represents an element which restricts to $(\gamma_5, 0)$ in $(H^*(2^4)^{\mathcal{A}_6})^2$. Also, $\gamma = \{\lambda_2 b_4 + x_3\lambda_3\}$ must represent an element which restricts to $(\bar{w}_3^2, 0)$ in $(H^*(2^4)^{\mathcal{A}_6})^2$. We have that $a_7$ is represented by $\lambda_2 x_5$ which implies the existence of the cup product $x_5 \cup a_2 = a_7$. Finally, $b_4^2$ represents $d_8$ and $d_6^2$ represents $d_{12}$.

We have not checked to see if there are higher differentials in the Eilenberg-Moore spectral sequence that we are using, but by inspection of the table above it is clear that through dimension 8 all the classes remaining at $E_3$ are infinite cycles. However, we will see now that the results to this point suffice to completely determine $H^*(PSU_4(3))$.

The first thing to note is that the two classes in dimension 5 have $Sq^4$'s on them which construct classes which restrict to $(\gamma_9, 0)$ and $(\gamma_9, \gamma_9)$ respectively. $a_7$ and $a_{11}$ must exist in the radical because $d_8$ and $d_{12}$ must be in the images of some higher order Bocksteins. The ring structure now implies that $\mathbb{F}_2[d_8, d_{12}](a_2, a_7, a_{11}) \subset \text{rad}$.

Next we write down all the elements in dimensions 15 and 16. In 16 there are 5 elements, $d_8^2, d_8\bar{w}_3\gamma_5, d_8\bar{w}_3\tau_5, \bar{w}_3^2\gamma_5^2$, and $\bar{w}_3^3\gamma_5\tau_5$, the second and third with non-trivial $Sq^1$'s and the fourth and fifth in the image of $Sq^1$ ($\gamma_5^3$ and $\gamma_5^2\tau_5$ respectively). The first image is in the Bockstein image from $d_8 a_7$. In 15 there are either two elements corresponding to $b_{15}$ (restricting respectively to $(b_{15}, 0)$ and $(0, b_{15})$ or there is just one which restricts to $(b_{15}, b_{15})$. The remaining elements which we are sure of are

$$\{d_{12}\bar{w}_3, \bar{w}_3^5, \bar{w}_3^3\tau_6, \gamma_5^3, \gamma_5^2\tau_5, \gamma_9\bar{w}_3^2, \tau_9\bar{w}_3^2\}.$$

In this list the third and fourth have $Sq^1$ images in 16, while the second is $Sq^1(\bar{w}_3^3\gamma_5)$. Also, the fifth is $Sq^1(\gamma_9\gamma_5)$ while the last is $Sq^1(\gamma_9\tau_5)$. Finally, the first is in the Bockstein image from $d_{12}a_2$.

Similarly, in 14 there are the four elements listed above, $d_8\bar{w}_3^2$, $d_8\tau_6$ and (at most) $a_{14}$. But the two elements involving $d_8$ are in the image of $Sq^1$ it follows that $a_{14}$ must be present to have a Bockstein to $b_{15}$, and there can only be the one element $b_{15}$ in dimension 15.

At this point it is probably convenient to write down the Poincaré series for the quotient by rad, at least as it is determined so far. So let us assume that the quotient ring of doubled elements is $\mathbb{F}_2[d_8, d_{12}](1, \bar{w}_3, b_{15}, \bar{w}_3 b_{15})$. Thus, $H^*(PSU_4(3))/\mathrm{rad}$ will have Poincaré series

$$\frac{2(1+x^9)(1+x^{15})}{(1-x^3)(1-x^5)(1-x^8)(1-x^(12))} - \frac{(1+x^3)(1+x^{15})}{(1-x^8)(1-x^{12})}$$
$$= \frac{(1+x^{15})(1-x^5-x^6+2x^9+x^{11})}{(1-x^3)(1-x^5)(1-x^8)(1-x^{11})}$$

which expands into a Taylor series of the form

$$1+x^3+2\,x^5+2\,x^6+3\,x^8+4\,x^9+2\,x^{10}+3\,x^{11}+5\,x^{12}+4\,x^{13}+6\,x^{14}+8\,x^{15}+5\,x^{16}$$
$$+10\,x^{17}+11\,x^{18}+7\,x^{19}+15\,x^{20}+16\,x^{21}+12\,x^{22}+18\,x^{23}+22\,x^{24}+O(x^{25}).$$

Actually, at this point there is one more possible doubled element, $\bar{w}_3 b_{15}$ which we assumed was doubled in writing down the Poincaré series above, but which we haven't yet actually proved is doubled. We do that now.

From the Poincaré series, in dimension 19 a total list of the elements in the quotient by the radical is

$$\{d_8^2\bar{w}_3, d_8\bar{w}_3^2\gamma_5, d_8\bar{w}_3^2\tau_5, \bar{w}_3^3\gamma_5^2, \bar{w}_3^3\gamma_5\tau_5, \gamma_9\gamma_5^2, \gamma_9'\gamma_5^2\}$$

and we see that two of them are in the image of $Sq^1$, $\gamma_9\bar{w}_3^3 \mapsto \bar{w}_3^3\gamma_5^2$, $\gamma_9'\bar{w}_3^3 \mapsto \bar{w}_3^3\gamma_5\tau_5$, $d_8^2\bar{w}_3$ is in the image of a higher Bockstein from $d_8^2 a_2$, and the remaining elments have non-trivial $Sq^1$'s.

In dimension 18 the 11 elements are

$$\{d_{12}\bar{w}_3^2, d_{12}\gamma_6, b_{15}\bar{w}_3, \gamma_9\bar{w}_3^3, \gamma_9\bar{w}_3\gamma_6, \gamma_5^3\bar{w}_3, \gamma_5^2\bar{w}_3\tau_5, d_8\gamma_5^2, d_8\gamma_5\tau_5, \bar{w}_3^6, \bar{w}_3^4\gamma_6\},$$

and two elements are needed for $Sq^1$ to dimension 19. This leaves nine elments (actually, $Sq^1(\gamma_5^3\bar{w}_3) = Sq^1(\bar{w}_3^3\gamma_9)$, and $Sq^1(\bar{w}_3\gamma_9\gamma_6) = Sq^1(\bar{w}_3\gamma_5\gamma_9')$, so the 9 elements include 7 monomials and the two sums of the elements above with equal $Sq^1$'s which together form a basis).

The 10 elements in dimension 17 are

$$\{d_{12}\gamma_5, d_{12}\tau_5, d_8\gamma_9, d_8\gamma_9', \bar{w}_3^4\gamma_5, \bar{w}_3^4\tau_5, \gamma_9\bar{w}_3\tau_5, \gamma_9\bar{w}_3\tau_5, d_8\bar{w}_3^3, d_8\bar{w}_3\gamma_6\}$$

and the last two elements are in the image of $Sq^1$ while the first eight map under $Sq^1$ to each of the elements in the basis for $\mathrm{Ker}(Sq^1)$ except $b_{15}\bar{w}_3$. It follows that only $b_{15}\bar{w}_3$ (and possibly an element which restricts to $(b_{15}\bar{w}_3, 0)$) are left unaccounted for in the Bockstein spectral sequence.

On the other hand there is one and only one element in the radical, $d_8a_{11}+d_{12}a_7$, which is in the kernel of the higher Bocksteins taking $a_7$ to $d_8$ and $a_{11}$ to $d_{12}$. It follows that $b_{15}\bar{w}_3$ must be doubled, must have a higher Bockstein to $d_8a_{11}+d_{12}a_7$ and there can be no element in $H^*(PSU_4(3))$ which restricts to $(b_{15}\bar{w}_3, 0)$.

The arguments above complete the determination of $H^*(PSU_4(3))$ and we have

THEOREM 3.7: *$rad(H^*(PSU_4(3))) \cong rad(H^*(L_3(4))) \cap H^*(S)^{\mathbb{Z}/2}$ which is the module $\mathbb{F}_2[d_8, d_{12}](a_2, a_7, a_{11}, a_{14})$ and*

$$H^*(PSU_4(3))/rad \xrightarrow{res^*} H^*(2^4_I)^{\mathcal{A}_6} \oplus H^*(2^4_{II})^{\mathcal{A}_6}$$

*is injective where the double image classes form the quotient algebra*

$$\mathbb{F}_2[d_8, d_{12}](1, \bar{w}_3, b_{15}, b_{15}\bar{w}_3).$$

*In particular, there is a class $\tau_5 \in H^5(PSU_4(3))$ with $res^*(\tau_5) = (\gamma_5, 0)$. Consequently, the same is true for $Sq^1(\tau_5) = \gamma_6$ which restricts to $(Sq^1(\gamma_5), 0) = (w_3^2, 0)$, and $Sq^4(\tau_5) = \gamma'_9$ which restricts to $(\gamma_9, 0)$.* ∎

## §4 The Bockstein spectral sequence for $PSU_4(3)$

In order to clarify the relations between the various elements in $H^*(PSU_4(3))$ as described in 3.7 we now analyze the Bockstein spectral sequence for $H^*(PSU_4(3))$.

To begin we write out the ideal, Ker$(res^*)$ where $res^*$ is restriction to $H^*(2^4_I)$ as follows:

$$\begin{aligned} I = \mathbb{F}_2[d_8, d_{12}](a_2, a_7, a_{11}, a_{14}) &\oplus x_5\mathbb{F}_2[\bar{w}_3, \gamma_5](1, \gamma_9, b_{15}, \gamma_9 b_{15}) \\ &\oplus \gamma_6\mathbb{F}_2[\bar{w}_3](1, \gamma_9, b_{15}, \gamma_9 b_{15}) \\ &\oplus \gamma'_9(1, \bar{w}_3, b_{15}, b_{15}\bar{w}_3)) \end{aligned}$$

The relations are $x_5\gamma'_9 = x_5\gamma_9 = \gamma_5\gamma'_9$, and $\gamma_6\gamma_5 = \bar{w}_3^2 x_5$, $\gamma_6\gamma'_9 = \bar{w}_3^2\gamma'_9$, $\gamma_6 x_5 = \gamma_6\gamma_5 = \bar{w}_3^2 x_5$. The first differential in the Bockstein spectral sequence is $Sq^1$, and here, $Sq^1(\gamma_5) = \bar{w}_3^2$, $Sq^1(x_5) = \gamma_6$, $Sq^1(\gamma^9) = \gamma_5^2$, and $Sq^1(\gamma'_9) = x_5\gamma_5$. Consequently we have

LEMMA 4.1: *In the ideal above, the $E_2$-term of the Bockstein spectral sequence has the form $\mathbb{F}_2[d_8, d_{12}](a_2, a_7, a_{11}, a_{14})$.*

PROOF: We have that the kernel of $Sq^1$ on $x_5\mathbb{F}_2[\bar{w}_3, \gamma_5, d_8, d_{12}](1, \gamma_9, b_{15}, \gamma_9 b_{15})$ is exactly $x_5\gamma_5\mathbb{F}_2[\bar{w}_3, \gamma_5^2, d_8, d_{12}](1, b_{15})$, while $x_5\mathbb{F}_2[\bar{w}_3, d_8, d_{12}](1, \gamma_9, b_{15}, \gamma_9 b_{15})$ has image

$$\gamma_6\mathbb{F}_2[\bar{w}_3, d_8, d_{12}](1, \gamma_9, b_{15}, \gamma_9 b_{15})$$

(summed with $x_5\mathbb{F}_2[\bar{w}_3, \gamma_5^2, d_8, d_{12}](1, b_{15})$ from the differential on $\gamma_9$, but we ignore that here, and simply regard the differential as removing the terms involving $\gamma_6$ from the list). Now, in the full submodule involving multiples of $x_5\gamma_5$,

$$x_5\gamma_5\mathbb{F}_2[\bar{w}_3, \gamma_5, d_8, d_{12}](1, \gamma_9, b_{15}, \gamma_9 b_{15})$$

we have that $x_5\gamma_5\gamma_9 \mapsto x_5\gamma_5(\gamma_5^2)$, $x_5\gamma_5(\gamma_5) \mapsto x_5\gamma_5\bar{w}_3^2$, in the kernel of $Sq^1$. Thus, after factoring the kernel of $Sq^1$ by these images, we have $x_5\gamma_5\mathbb{F}_2[d_8, d_{12}](1, \bar{w}_3, b_{15}, \bar{w}_3 b_{15})$ left. On the other hand, the image under $Sq^1$ of the part involving $\gamma_9'$ is exactly this, so we get cancellation of all these terms, leaving just the terms asserted in the lemma. ∎

We also have, as a consequence of the calculations above

LEMMA 4.2: *In the quotient by the ideal above,* $\mathbb{F}_2[\bar{w}_3, \gamma_5, d_8, d_{12}](1, \gamma_9, b_{15}, \gamma_9 b_{15})$ *the contribution to the* $E_2$*-term of the Bockstein spectral sequence is* $\mathbb{F}_2[d_8, d_{12}](1, \bar{w}_3, b_{15}, \bar{w}_3 b_{15})$.

∎

Thus, we conclude that the entire $E_2$-term of the Bockstein spectral sequence has the form

4.3 $$E_2 = \mathbb{F}_2[d_8, d_{12}](1, a_2, \bar{w}_3, a_7, a_{11}, a_{14}, b_{15}, \bar{w}_3 b_{15})$$

We know that $d_2(a_2) = \bar{w}_3$. Also, an early differential on $a_{14}$ must be $b_{15}$. At this point, the spectral sequence will have the form

$$\mathbb{F}_2[d_8, d_{12}](1, a_7, a_{11}, \bar{w}_3 b_{15}).$$

Then, the remaining differentials must have the form $d(a_7) = d_8$, $d(a_{11}) = d_{12}$, and $d(\bar{w}_3 b_{15}) = d_8 a_{11} + d_{12} a_7$, though we have no idea of the order in which they occur.

REMARK 4.4: As a consequence of the last calculation we have the the higher 2-torsion in $H^*(PSU_4(3); \mathbb{Z})$ when tensored with $\mathbb{F}_2$ has the form $\mathbb{F}_2[d_8, d_{12}](1, \bar{w}_3, b_{15}, d_8 a_{11} + d_{12} a_7)$.

## §5 The Determination of $H^*(McL)$

The calculation of $H^*(McL))$ will be divided into 2 steps. Namely, we will determine its non–nilpotent part and its radical separately. To start, we have the exact sequence

$$0 \to \mathrm{rad}(H^*(McL)) \longrightarrow H^*(McL) \to H^*(2_I^4)^{\mathcal{A}_7} \oplus H^*(2_{II}^4)^{\mathcal{A}_7}. \qquad (5.1)$$

As pointed out in 2.4, the image of $res_{E_i}^{McL}$ will be the full ring of invariants in each case. From 2.6,

$$H^*(2^4)^{\mathcal{A}_7} = \mathbb{F}_2[d_8, d_{12}, d_{14}, d_{15}](1, x_{18}, x_{20}, x_{21}, x_{24}, x_{25}, x_{45})$$

so, in particular, $H^*(McL) \cong \mathbb{F}_2[d_8, d_{12}] \oplus \mathrm{rad}$ through dimension 13, and contains $\mathbb{F}_2[d_8, d_{12}]$ as a sub-polynomial algebra. It follows for Bockstein reasons that $a_7$ and $a_{11}$ must both be contained in the radical, so $\mathbb{F}_2[d_8, d_{12}](a_7, a_{11}) \subset \mathrm{rad}$. Hence, from 2.2 we have

LEMMA 5.1: *The radical of* $H^*(McL)$ *is exactly* $\mathbb{F}_2[d_8, d_{12}](a_7, a_{11})$. ∎

We now turn to the double image classes. From the inclusion $PSU_4(3) \subset McL$ we have the commutative diagram of restriction maps,

$$\begin{array}{ccc} H^*(McL) & \xrightarrow{\coprod res^*} & H^*(2^4_I)^{\mathcal{A}_7} \oplus H^*(2^4_{II})^{\mathcal{A}_7} \\ \Big\downarrow{\scriptstyle res^*} & & \Big\downarrow{\scriptstyle in\oplus in} \\ H^*(PSU_4(3)) & \xrightarrow{\coprod res^*} & H^*(2^4_I) \oplus H^*(2^4_{II}) \end{array}$$

where the bottom restriction map was discussed in §3. As before, to make this meaningful, we must calculate the double image classes i.e. pairs $(res^{McL}_{2^4_I}(x), res^{McL}_{2^4_{II}}(x))$ where both components are non–trivial. To do this we use the double coset decomposition of [JW]:

5.2 $$McL \;=\; PSU_4(3) \sqcup PSU_4(3)xPSU_4(3) \sqcup PSU_4(3)v_1PSU_4(3)$$

where $x$ has order 2 and the order of $v_1$ is four. The two intersections are $L_3(4) = xPSU_4(3)x \cap PSU_4(3)$ and $3^4{:}\mathcal{A}_6$ where the action of $v_1$ leaves invariant the $D_8$-Sylow subgroup, but acts on it as the non-trivial outer automorphism. This only produces a stability condition on the restrictions to $2^2$'s in $McL$. The double coset involving $x$ produces the copy of $L_3(4){:}2_2 \subset McL$, and additionally gives the extra automorphisms to build the $\mathcal{A}_7$ action on the two $2^4$'s from the $\mathcal{A}_6$-action in $PSU_4(3)$. But these are then the *only* stability requirements which are placed on the elements of $H^*(PSU_4(3))$ in order that they be in the image from $H^*(McL)$.

Corollary 5.3: *Let $\alpha \in H^*(PSU_4(3))$, then $\alpha$ is contained in the image of $res^*$ from $H^*(McL)$ if $\alpha$ restricts to $(\lambda_1, \lambda_2) \in H^*(2^4_I)^{\mathcal{A}_7} \oplus H^*(2^4_{II})^{\mathcal{A}_7}$ and $\lambda_i$, $i = 1, 2$ restricts to zero in $H^*(2^2)$ for some $2^2 \subset 2^4$ provided the $\lambda_i$ are in a dimension where the radical of $H^*(PSU_4(3))$ is equal to the radical of $H^*(McL)$.*

Proof: The action of $\mathcal{A}_7$ on $2^4$ is transitive on the 35 distinct subgroups of the form $2^2 \subset 2^4$. Thus, for any $\alpha \in H^*(2^4)^{\mathcal{A}_7}$, in order to verify that $\alpha$ restricts to zero in $H^*(2^2)$ for each $2^2 \subset 2^4$ it suffices to check on a single subgroup. Moreover, since the radical of $H^*(PSU_4(3))$ is equal to the radical of $H^*(McL)$ in this dimension it follows that $\alpha$ is completely determined by its restrictions to $H^*(2^4_I)$ and $H^*(2^4_{II})$. Consequently, it is stable for the double coset decomposition in 5.2. ∎

Remark: The only dimensions where there are elements in $\mathrm{rad}(H^*(PSU_4(3)))$ which are not in $\mathrm{rad}(H^*(McL))$ are of the form $10, 14, 18, 22, 26, 30$ in the range we will need to study.

We are now ready to determine $H^*(McL)$.

Theorem 5.4: *Let $\mathcal{M}$ be the algebra obtained by identifying the two copies of the subalgebra $\mathbb{F}_2[d_8, d_{12}](1, x_{18})$ in the direct sum of two copies of the algebra of $\mathcal{A}_7$-invariants*

$$\mathbb{F}_2[d_8, d_{12}, d_{14}, d_{15}](1, x_{18}, x_{20}, x_{21}, x_{24}, x_{25}, x_{27}, x_{45}).$$

*Then $H^*(McL)$ is described as an extension*

$$0 \to \mathbb{F}_2[d_8, d_{12}](a_7, a_{11}) \to H^*(McL) \to \mathcal{M} \to 0.$$

In other words, we have a short exact sequence

$$0 \to \mathcal{M} \to H^*(2^4_I)^{\mathcal{A}_7} \oplus H^*(2^4_{II})^{\mathcal{A}_7} \to \mathbb{F}_2[d_8, d_{12}](1, x_{18}) \to 0$$

where $\mathcal{M}$ is isomorphic to the algebra $H^*(McL)/\mathrm{rad}$.

Proof: Note that in dimensions 20, 21, 24 and 25 there are no elements in the radical of $H^*(PSU_4(3))$. Also, in [AM4] we write down representatives for the four elements $x_{20}$, $x_{21}$, $x_{24}$ and $x_{25}$ in 2.4 in the form of symmetric sums. Note that each symmetric sum appearing in the expressions for these elements has the form $Sx_1^i x_2^j x_3^k$ or $Sx_1^i x_2^j x_3^k x_4^l$ with $i, j, k, l$ all positive. Consequently, restricting each of these elements to the $2^2$ spanned by elements dual to $x_1$ and $x_2$ does give zero, and $(x_{20}, 0)$, $(x_{21}, 0)$, $(x_{24}, 0)$ and $(x_{25}, 0)$ are all in the image from $H^*(McL)$ if they were not represented by doubled elements in $H^*(PSU_4(3))$. But from 3.7 there are no double image elements in these dimensions except for $d_8 d_{12}$ which has already been accounted for. It follows that $(x_{20}, 0)$, $(x_{21}, 0)$, $(x_{24}, 0)$ and $(x_{25}, 0)$ are all in the restriction image from classes in $H^*(McL)$.

We now show that the classes $(d_{14}, 0)$ and $(d_{15}, 0)$ are in the restriction image from classes in $H^*(McL)$. From [AM4, 3.14] we have $d_{14} = \gamma_5\gamma_9 + \bar{w}_3^2 d_8 + \bar{w}_3^3 \gamma_5$ from which we see that there is also the element $\bar{d}_{14} = \gamma_5\gamma_9' + \bar{w}(\gamma_6 d_8 + \bar{\gamma}_5) \in H^*(PSU_4(3))$ which restricts to $(d_{14}, 0)$ in the two invariant subrings as well as restricting trivially to $H^*(2^2)$. There is, however, the element $d_{12}a_2$ in the radical in this dimension so we do have to be a bit careful about asserting that $\bar{d}_{14}$ actually lies in $H^*(McL)$. However, in dimension 15, the only element in the radical is $d_8 a_7$ which, as we have seen, is in $H^*(McL)$, and a similar argument constructs $\bar{d}_{15}$ which restricts to $(d_{15}, 0)$, and which therefore does represent an element in $H^{15}(McL)$. Consequently, $\bar{d}_{14}$ must also be present in order to have the Bocksteins work out correctly. In particular, $Sq^1(d_{14}) = d_{15}$ and $Sq^1(\bar{d}_{14}) = \bar{d}_{15}$.

From the explicit descriptions of the classes $x_i$ in [AM4] we have $Sq^1(x_{20}) = x_{21}$, $Sq^1(x_{24}) = x_{25}$ in $H^*(2^4)^{\mathcal{A}_7}$. This gives us sufficient information to understand the first differential $Sq^1$ in the Bockstein spectral sequence for $H^*(McL)$, at least through dimension 27. We see that, including the class $x_{27}$ and the (at this stage only possible) class $\bar{x}_{27}$ which restricts to $(x_{27}, 0)$, the $E_2$-term has the form

$$\mathbb{F}_2[d_8, d_{12}, d_{14}^2](1, x_{18}, x_{27}, \bar{x}_{27}, d_{14}\bar{d}_{14}, x_{18}d_{14}\bar{d}_{14}) \oplus \mathbb{F}_2[d_8, d_{12}](a_7, a_{11})$$

at least assuming that $Sq^1(x_{27}) = 0$. But from this it is clear that $\bar{x}_{27}$ must be present to kill $d_{14}\bar{d}_{14}$, and that $x_{18}$ must be doubled – there cannot be a class $\bar{x}_{18}$ restricting to $(x_{18}, 0)$ in $H^*(McL)$ since there are no classes in either dimensions 17 or 19 to account for it in the Bockstein spectral sequence. Finally, the role of the class $x_{45}$ must be to hit $x_{18}d_{14}^2$ but since the class $x_{18}d_{14}\bar{d}_{14}$ is also present in dimension 46 it follows that there is a class $\bar{x}_{45} \in H^*(McL)$ which restricts to $(x_{45}, 0)$.

This completes the construction of the required classes and complete the proof. ∎

COROLLARY 5.5: *The Poincaré series for $H^*(McL)$ is given by*

$$p(t) = \frac{2(1+t^{18}+t^{20}+t^{21}+t^{24}+t^{25}+t^{27}+t^{45})}{(1-t^8)(1-t^{12})(1-t^{14})(1-t^{15})} + \frac{(t^7-1)(1-t^{11})}{(1-t^8)(1-t^{12})}.$$

∎

Note in particular that we have that $H^i(McL)=0$ for $1 \le i \le 6$ and that $H^7(McL) \cong H^8(McL) = \mathbb{F}_2$.

## §6 The Extension $3^{\cdot}McL.2$

According to the Atlas [Co], there is a subgroup $N$ of the form $3^{\cdot}McL.2$ which is maximal and of odd index in the sporadic group $Ly$. The new element of order 2 exchanges $2^4_I$ and $2^4_{II}$, and acts trivially on $\mathrm{rad}(H^*(McL))$. From this we deduce first that $2^4_I \subset Syl_2(Ly) \subset Ly$ is weakly closed and that the $E_2$ term of the Lyndon–Hochschild–Serre spectral sequence associated to the extension has the form

$$\mathbb{F}_2[d_8,d_{12}](a_7,a_{11}) \oplus H^*((\mathbb{Z}/2)^4)^{\mathcal{A}_7} \oplus \mathbb{F}_2[d_8,d_{12},e](1,a_7,a_{11},x_{18})e.$$

Indeed, only the doubled elements can multiply non–trivially with the one dimensional generator $e$.

PROPOSITION 6.1: *The spectral sequence above collapses at the $E_2$–stage, i.e. all differentials are trivial.*

PROOF: For Bockstein reasons (connecting $a_7$ to $d_8$, $a_{11}$ to $d_{12}$) and the fact, due to weak closure, that the invariants $H^*(2^4_i)^{\mathcal{A}_7}$ are in the image of the restriction map, it suffices to show that $d_8$, $d_{12}$ are permanent cocycles. From the Atlas [Co], page 175, there exists a real representation $\chi_4$ for $Ly$ and hence for $3^{\cdot}McL.2$ which restricted to $(\mathbb{Z}/2)^4$ has the form $110\epsilon \oplus 2849\chi_{reg}$, where $\epsilon$ is the trivial representation, and $\chi_{reg}$ is the regular representation. Hence $res(w_8(\chi_4)) = d_8$, and $res(w_{12}(\chi_4)) = d_{12}$, showing that they are permanent cocycles in the spectral sequence. ∎

In the Atlas [Co], the double coset decomposition of $G = Ly$ with respect to $N = 3^{\cdot}McL.2$ is explicitly described. In fact we have

$$G = N \cup Nx_1N \cup Nx_2N \cup Nx_3N \cup Nx_4N$$

with known intersections: $3^{2+4}.4S_5$, $2S_7$, $4S_6$ and $5^{1+2}:S_3$ respectively. This will be used to calculate $H^*(Ly)$ in a sequel.

## References

[AM1] Adem, A. and Milgram, R. J. *Cohomology of Finite Groups*, Springer–Verlag Grundlehren 309 (1994).

[AM2] Adem, A. and Milgram, R. J., "The Cohomology of the Mathieu Group $M_{22}$," *Topology*, to appear.

[AM3] Adem, A. and Milgram, R. J., "$\mathcal{A}_5$–Invariants, the Cohomology of $L_3(4)$ and Related Extensions," *Proc. London Math. Society* (3) 66 (1993), 187–224.

[AM4] Adem, A. and Milgram, R. J., "Invariants and Cohomology of Groups," *Bol. Soc. Mat. Mex.*, to appear.

[Co] Conway, J. et al, *Atlas of Finite Groups*, Oxford University Press 1985.

[FP] Fiedorowicz, S. and Priddy, S., *Homology of Classical Groups over Finite Fields and their Associated Infinite Loop Spaces* Springer–Verlag Lecture Notes in Mathematics 674 (1978).

[J] Janko, Z., "A Characterization of the Mathieu Simple Groups," *Journal of Algebra* 9 (1968) 20–41.

[JW] Janko, Z. and Wong, W., "A Characterization of McLaughlin's Simple Group," *Journal of Algebra* 20 (1972) 203–225.

[McL] McLaughlin, J., "A simple group of order 898,128,000," *The Theorey of Finite Groups*, (Brauer and Sah, eds.), Benjamin, 1969, pps. 109-111.

[M] Milgram, R.J., "The cohomology of the Mathieu group $M_{23}$," (Preprint, Stanford, 1993).

[RSY] A. Ryba, S. Smith and S. Yoshiara, "Some projective modules determined by sporadic geometries," J. Algebra 129 (1990), 279–311.